\documentclass{article}
\usepackage[utf8]{inputenc}
\usepackage[frozencache,cachedir=.]{minted}
\usepackage{graphicx}
\usepackage{hyperref}
\usepackage{amsmath}
\usepackage{amssymb}
\usepackage{fullpage}

\usepackage[
backend=biber,
style=ieee,
sorting=none
]{biblatex}

\newcommand{\toolbox}{TALOS}
\newcommand{\toolboxfullname}{Toolbox for Analysis and Large-scale Optimization of Spacecraft}
\newcommand{\model}{\texttt{Model}}

\title{\toolbox: A toolbox for spacecraft conceptual design}
\author{Victor Gandarillas and John T. Hwang}
\date{}

\begin{document}

\maketitle

\begin{abstract}
We present the \toolboxfullname{} (\toolbox{}), a framework designed for applying large-scale multidisciplinary design optimization (MDO) to spacecraft design problems.
The framework is built using the Computational System Design Language (CSDL), with abstractions for users to describe systems at a high level.
CSDL is a compiled, embedded domain-specific language that fully automates derivative computation using the adjoint method.
CSDL provides a unified interface for defining MDO problems, separating model definition from low-level program implementation details.
\toolbox{} provides discipline models for spacecraft mission designers to perform analyses, optimizations, and trade studies early in the design process.
\toolbox{} also provides interfaces for users to provide high-level system descriptions without the need to use CSDL directly, which simplifies the exploration of different spacecraft configurations.
We describe the interfaces in \toolbox{} available to users and run analyses on selected spacecraft subsystem disciplines to demonstrate the current capabilities of \toolbox{}.
\end{abstract}

\section{Introduction}\label{sec:intro}

The goal of the \toolboxfullname{} (\toolbox{}) is to facilitate modeling, analysis, and multidisciplinary design optimization (MDO) of spacecraft vehicle and mission concepts.
The emphasis of \toolbox{} is currently on small satellite (SmallSat) and CubeSat design, but the library is extensible enough to accommodate general spacecraft design problems.

SmallSats and CubeSats are classes of miniature satellites that are becoming
increasingly popular for scientific and commercial missions due to their
low cost and versatility.
SmallSats and CubeSats have been used for a variety of missions, including Earth
observation and scientific research.
SmallSats and CubeSats offer a relatively affordable alternative to
traditional large-scale missions and allow for increased frequency of
launches.
In recent years, SmallSats and CubeSats have seen a significant
growth in mission diversity, including Earth observation, planetary
exploration, space technology development, and scientific research.

A diversity of mission requirements has led to a diversity of SmallSat designs as well as CubeSat designs.
The SmallSat X-ray Quantum Calorimeter Satellite
(XQCSat)~\cite{mccammon2019xqcsat, kaaret2019xqcsat} is designed to
perform X-ray spectroscopy.
CubeSats such as the Radiometer Assessment using Vertically Aligned
Nanotubes (RAVAN)~\cite{swartz2019radiometer}, designed to measure the
Earth's outgoing radiation, and the Lunar
Flashlight~\cite{cohen2020lunar}, designed to search for ice on
the Moon's surface, also showcase the diversity of satellite designs.
The Lunar Flashlight mission demonstrates how CubeSats can be used for
planetary exploration missions, which has traditionally been the domain
of larger, more expensive spacecraft.

One advantage of SmallSats and CubeSats is their small size, which
makes them easy to launch and allows for multiple spacecraft to be
launched at once.
Some satellite missions require multiple spacecraft to fly in
formation.
The TerraSAR-X and Tandem-X (TerraSAR-X add-on for Digital Elevation
Measurement)~\cite{krieger2007tandem} spacecraft use Synthetic Aperture
Radar (SAR) to produce high-resolution images of the Earth's surface and
generate a Digital Elevation Model (DEM) of the Earth's surface.
The TechSat-21~\cite{burns2000techsat} was a technology demonstration
mission focused on demonstrating new technologies for formation flight.
TanDEM-X and TechSat-21 are both remote sensing systems, but SmallSats
have filled the need for formation flight beyond Earth orbit as well.
The Double Asteroid Redirection Test (DART)~\cite{cheng2012dart,
rivkin2021double} was a technology demonstrator aimed at redirecting
near-Earth objects (NEOs) that successfully redirected an asteroid.

Virtual telescopes
represent another application of SmallSats and CubeSats, which leverage formation-flying
to observe celestial objects.
These platforms
offer a cost-effective alternative to traditional telescopes, which can
be expensive to build and maintain.
Virtual telescopes are formed from SmallSats, such as the
Miniature Distributed Occulter Telescope
(mDOT)~\cite{macintosh2019miniature}, Virtual Telescope for X-ray
Observations (VTXO)~\cite{rankin_vtxo_nodate}, and
Proba-3~\cite{llorente2013proba, galano2019proba} missions, or CubeSats,
such as the VIrtual Super-resolution Optics with Reconfigurable Swarms
(VISORS)~\cite{lightsey2022concept, gundamraj2021preliminary,
koenig2021formation}.

These missions highlight the growing importance of SmallSats and
CubeSats in remote sensing and space exploration. This trend is likely
to continue in the future, as SmallSats offer a flexible and
affordable platform for a variety of scientific and technological
objectives.
The complexity and diversity of SmallSat and CubeSat mission
profiles highlight the need for rapid model development for generating
new configurations for vehicle concepts and performing trade studies.
Systems Tool Kit (STK) and FreeFlyer are commercial software packages
that are widely used for spacecraft mission planning, design, and
analysis.
Although these software applications built for spacecraft design are capable of high-fidelity
modeling, mission designers could benefit from a  framework that relies
on optimization techniques to generate vehicle concepts and conduct
trade studies.

Multidisciplinary design optimization (MDO) is the use of numerical
optimization to generate designs for engineering systems that involve multiple subsystems or
disciplines~\cite{martins2013multidisciplinary}.
Many MDO problems have objective and constraint functions that are
continuous and differentiable in the design variables.
These problems are also known as nonlinear programs.
Gradient-based approaches to solving nonlinear programs scale better
than gradient-free approaches.
A disadvantage of gradient-based approaches is that they require
derivative information from the model of the system, which can be
difficult to implement for large-scale multidisciplinary systems, especially
when there is coupling across disciplines.
Computing exact derivatives for models with multiple disciplines relies
on computational architectures~\cite{hwang2018computational} and
software frameworks~\cite{gray2019openmdao} that facilitate
multidisciplinary model definition, evaluation, and derivative
computation.
The Computational System Design Language
(CSDL)~\cite{gandarillas2022methodology} targets solving MDO
problems and fully automates derivative computation using the adjoint
method.

\toolbox{} is written in CSDL with the goal of providing users with a simplified interface for rapidly generating multidisciplinary models to run analyses and trade studies as well as solving large-scale MDO problems with minimal programming effort.
The main components of the \toolbox{} library are predefined numerical models written using CSDL and interfaces for users to provide high-level system descriptions called specifications.
CSDL simplifies numerical modeling and generates a low-level implementation of a computational model.
As a result, mathematical modeling and program implementation (including any implementation of a numerical algorithms necessary to evaluate the model) are completely separate, allowing engineers to focus on a model of the physical system, and not on the implementation of a program that simulates the behavior of the modeled system.
The generated program also fully automates derivative computation using the adjoint method so that all models defined using CSDL can be used within an optimization framework using a gradient-based approach.
As a result, the same model definition can be used in both analysis and optimization contexts with no additional overhead from the user.

In addition to leveraging the simplicity of CSDL to define numerical models, \toolbox{} provides interfaces for users to provide high-level system descriptions called specifications without the need to describe mathematical relationships between physical quantities in the model.
These high-level descriptions offer a higher level of abstraction than what CSDL can offer, and automate the assembly of complex model hierarchies in CSDL.
This article introduces users to the various components available in \toolbox{} for designing a CubeSat mission.

\section{Software Design}

CSDL was selected as the language for \toolbox{} due its ability to fully automate derivative computation as well as the use of both the functional and object-oriented paradigms.
The functional paradigm allows users to define models using natural syntax that closely resembles mathematical notation, while the object-oriented paradigm allows users to compose models hierarchically.
These features provide users with the ability to apply gradient-based methods to optimization problems with minimal overhead.

\toolbox{} further simplifies model definition by offering a higher level of abstraction over spacecraft design than what is available by using CSDL directly.
\toolbox{} achieves superior usability, extensibility, and portability, summarized in the following subsections.

\subsection{Usability}

\toolbox{} enables developers with high level knowledge of space systems and limited programming experience to assemble a multidisciplinary model with ease.
Users specify a system configuration by defining a high level system description using the concept of a specification that \toolbox{} provides.
Specifications direct \toolbox{} to build model hierarchies from prepackaged models defined in CSDL without requiring the user to define models directly using CSDL.

All models in \toolbox{} are written in CSDL, which simplifies development for contributors to \toolbox{}.
Interfaces for models of individual disciplines are clearly defined so that users can integrate discipline models into a single system level model.
The discipline models are also designed so that \toolbox{} specifications can control how discipline models are composed hierarchically.

\subsection{Extensibility}

\toolbox{} is extensible in the sense that new models and high-level system descriptions (specifications) can be added easily.
Model definition files are grouped according to whether they define a mission configuration or a spacecraft discipline.
Mission configurations generally contain spacecraft disciplines, while disciplines representing spacecraft subsystems generally have flatter hierarchies or contain no other discipline models.
A new model can be added to \toolbox{} by adding the file where the model is defined to the appropriate directory.
Model specifications can also be added, contain other specifications, and are all located in a directory separate from the model disciplines.

\subsection{Portability}

\toolbox{} is written in the Computational System Design Language (CSDL), a compiled, embedded domain-specific language.
CSDL inherits the simple syntax of its host language, Python.
The CSDL compiler is a three-stage compiler, which constructs an intermediate representation of a program from a given language, and generates a program in another language.
As a result, models defined in CSDL can theoretically be compiled to any low-level language, according to implementation requirements.
Building \toolbox{} on top of a language with a simple syntax that can be compiled to any language removes any limitation to which platform can be used to perform analyses or solve design optimization problems.
Furthermore, \toolbox{} specifications are written in pure Python, which can run on any modern desktop computer.

\section{Implementation}

\toolbox{} makes use of CSDL to define numerical models.
CSDL separates numerical model definition from computational model implementation.
The CSDL compiler generates a low-level implementation of the program used to simulate the model, hiding details from the user such as memory management and pointer arithmetic.
Even the implementation of numerical algorithms for solving equations that do not have closed-form solutions is transparent to the user.
\toolbox{} further increases the level of abstraction by providing predefined models specific to spacecraft design.
Users of \toolbox{} integrate predefined models provided by the library, and contributors to \toolbox{} need not implement low-level implementations.
As a result, \toolbox{} developers operate at a very high level of abstraction, focusing on the physical aspects of the design problem rather than the computational aspects.
Furthermore, \toolbox{} users are free to focus on the conceptual aspects of the design problem, focusing on defining a configuration rather than modeling the physical system.
This lends \toolbox{} to serve as a tool for rapid prototyping, making \toolbox{} ideal for conducting trade studies for spacecraft swarm missions.

Each model in \toolbox{} is defined as a CSDL \model{} subclass.
CSDL uses the functional programming paradigm for users to define models using a natural syntax closer to mathematical notation, and the object-oriented programming paradigm to enable hierarchical composition of models.
As a result of defining models in CSDL, \toolbox{} enables users to easily build spacecraft system models by hierarchically composing predefined models.

The primary interface \toolbox{} users interact with is the specification concept.
The motivation for the specification concept originates from the fact that the degree to which CSDL \model{} subclass definitions can vary in complexity.
Slight modifications to a highly complex \model{} subclass can easily lead to an explosion in the number of \model{} subclass definitions.
To minimize the number of \model{} subclass definitions and maximize code reuse in a given library, \model{} subclass constructors accept parameters that modify a model definition on a per-instance basis.

Parameters are compile time constants that a user can set to change the definition of a \model{} object wherever a \model{} subclass is instantiated.
Parameters can be of any Python type, and individual \model{} subclasses can use the full power of the Python language to programmatically control compile time behavior, namely, model definition.
For example, if the size of a variable in a \model{} class depends on the number of time steps, the user can pass an integer value to a \model{} subclass constructor to specify the size of the variables whose sizes depend on the number of time steps.
One use case would be to specify the number of time steps used in numerical integration.
Another example of a parameter controlling compile-time behavior would be a boolean flag that specifies whether or not a subsystem should be modeled for a given configuration.

Since parameters are useful for specifying model definition on a per-instance basis, collections of parameters used frequently across models can define model specifications.
For this reason, \toolbox{} provides specifications in the form of Python classes that contain a collection of parameters and/or other specifications.
Specifications themselves are Python types, and can therefore be treated as CSDL parameters.
The advantage of specifications is that they provide a single object to pass as a \model{} constructor argument that can then be passed to models further down the model hierarchy, simplifying the interfaces for \model{} subclasses provided by \toolbox{}.

\subsection{Tour of the \toolbox{} Library}\label{sec:tour}

This section presents the organization of the \toolbox{} library.
First, we present an overview of the organization of the \toolbox{} codebase to provide a sufficiently detailed view for a new user to easily find components appropriate for defining a particular  mission configuration.
Second, we present examples of analyses using some of the currently available models representing various spacecraft disciplines.
As \toolbox{} is under ongoing development, only some of the currently available models are used in the presentation of example analyses.
\toolbox{} is available for download at \href{https://github.com/lsdolab/talos}{https://github.com/lsdolab/talos}.

The library is organized hierarchically according to the type of object used to define a system model.
Multiple types of objects are available for defining a system model for a spacecraft mission configuration.
The core set of object types are the \model{} subclasses that represent the various subsystems relevant to space mission analysis and design.
Each \model{} subclass inherits from the \model{} class defined in the Computational System Design Language (CSDL).
\toolbox{} provides users with interfaces to modify model specifications through the use of CSDL parameters.
CSDL parameters allow model developers to modify model definitions when \model{} subclasses are instantiated as opposed to creating a new \model{} subclass each time a change in the model definition is necessary.
Parameters are allowed to be of any Python type, and \model{} subclasses can programmatically change their definition at compile time according to the parameters provided at construction.
This increases the reusability and maintainability of discipline model definitions provided by the \toolbox{} library.
To further simplify model definition and reuse, \toolbox{} also provides parameters implemented as data structures that contain collections of parameters.
These types of parameters are referred to as specifications in \toolbox{}.
Object types and classes are organized within \toolbox{} in the following directories:

\begin{itemize}
    \item \texttt{configurations/}
    \item \texttt{disciplines/}
    \item \texttt{examples/}
    \item \texttt{specifications/}
    \item \texttt{utils/}
\end{itemize}

Each directory contains type or class definitions according to their use case.
The \texttt{disciplines/} directory contains \model{} subclasses defining models of specialized disciplines.
Each discipline represents a subsystem on a spacecraft, or one aspect of the spacecraft environment.
For example, the \texttt{disciplines/} directory includes models of a spacecraft orbit, attitude dynamics, and propulsion system.
The \texttt{configurations/} directory contains models composed of other more specialized disciplines.
For example, the \texttt{configurations/} directory includes models of a CubeSat (including the relevant subsystems) and a virtual telescope (which also includes the CubeSat configuration as a subdiscipline).

The difference between the \texttt{disciplines/} directory and the \texttt{configurations/} directory is that the \texttt{disciplines/} directory is meant to provide components to assemble a configuration according to the needs of any given mission, while the \texttt{configurations/} directory is meant to provide a model of a multidisciplinary system suitable for a particular application.
The \texttt{configurations/} directory contains more complex models composed of models defined in the \texttt{disciplines/} directory.
The models in the \texttt{disciplines/} directory provide users with greater flexibility and control when building a system model, while the \texttt{configurations/} directory provides users with the ability to more easily and quickly reuse frequently used model hierarchies.

The \texttt{specifications/} directory contains data structures called specifications containing collections of parameters that fully specify a particular model.
Specifications are especially useful for passing many parameters from models at a higher level within the model hierarchy to models further down the model hierarchy in a clean way.
For example, the \texttt{specifications/} directory includes data structures containing parameters used to specify a ground station, a CubeSat, and a virtual telescope.
The details of how each discipline is defined and how the disciplines are hierarchically composed are hidden from the \toolbox{} user.
Each parameter in a specification controls certain parts of the model definitions such as the number of time steps used to numerically integrate differential equations or whether or not certain subsystems are modeled within a configuration.

The \texttt{utils/} directory contains functions not part of the CSDL Standard Library that facilitate model definitions.
For example, the \texttt{utils/} directory includes a ray tracing algorithm for generating data to train a surrogate model of the illumination of solar arrays as a function of spacecraft orientation.

\subsection{Specifications Provided by \toolbox{}}

\toolbox{} provides discipline models for spacecraft subsystems from which users can assemble system level models.
Many mission and vehicle concepts share discipline models, such as those described in Section~\ref{sec:disciplines}.
Users can thus benefit from repeating patterns of model composition when creating each new mission or vehicle concept.
\toolbox{} provides the following components called specifications that require users to provide only parameter values unique to a configuration, abstracting the composition of models:

\begin{itemize}
\item \texttt{GroundStationSpec}
\item \texttt{CubeSatSpec}
\item \texttt{VirtualTelescopeSpec}
\end{itemize}

The \texttt{GroundStationSpec} class contains parameters used to define a ground station to model communication between a spacecraft and the Earth.
\toolbox{} specifications can not only contain parameters, but also other specifications.
For example, if a user decides to model communication between a CubeSat and multiple ground stations located across the surface of the Earth, the \texttt{CubeSatSpec} can accept a dictionary of \texttt{GroundStationSpec} objects to indicate which ground stations are in communication with the CubeSat.
The \texttt{CubeSatSpec} class constructor arguments are the only details the user needs to provide in order to fully describe a CubeSat configuration in \toolbox{}.
Users do not need to redefine subsystem models, assemble a system model for a CubeSat configuration using the CSDL API, or issue connections between variables in different models within the model hierarchy.
Likewise, the \texttt{VirtualTelescopeSpec} requires two \texttt{CubeSatSpec} objects to specify the configurations for each CubeSat that forms the virtual telescope configuration.
The relationships between the CubeSat positions and the constraints on the telescope formation are predefined, and the user only needs to set the constraint values.
All of these specifications enable the \toolbox{} user to specify a configuration rather than a numerical model.
The following subsections show example use cases of each of these specifications.

\subsubsection{Ground Station Specification}\label{sec:gs}

The ground station specification provides \toolbox{} users with an interface to specify parameters for CubeSat communication discipline models without the need to integrate ground station discipline models with the system model in CSDL.
Listing~\ref{lst:gs} shows an example of how various ground stations are specified in \toolbox{}.
These specifications can be reused across other specifications as shown in~\ref{sec:cs} and \ref{sec:vt}.

\begin{listing}
\begin{minted}[]{python}
groundstations = {
    'SanDiego': GroundStationSpec(
        name='SanDiego',
        lon=-117.2340,
        lat=32.8801,
        alt=0.4849,
    ),
    'Atlanta': GroundStationSpec(
        name='Atlanta',
        lon=-84.3963,
        lat=33.7756,
        alt=0.2969,
    ),
}
\end{minted}
\caption{Creating a dictionary of ground station specifications. Other specifications also use the \texttt{GroundStationSpec} instances to construct discipline models for each ground station and communication between each ground station and each spacecraft in the system model.}
\label{lst:gs}
\end{listing}

\subsubsection{CubeSat Specification}\label{sec:cs}

The CubeSat specification provides the interface for defining a CubeSat configuration. Listing~\ref{lst:cs} shows an example of a specification defining a CubeSat configuration.
The \texttt{orbit\_model} parameter specifies which type of orbit model is to be used.
Two options are available: \texttt{`absolute'}, and \texttt{`relative'}.
Selecting \texttt{`absolute'} directs CSDL to use a model of the CubeSat's orbit that is a function of the CubeSat's position relative to the center of the Earth.
This option is suitable for mission configurations that use a single CubeSat, or CubeSats that do not need to maintain formation with tight tolerances.
Selecting \texttt{`relative'} directs CSDL to use a model of the CubeSat's orbit that is a function of the CubeSat's position relative to a reference orbit that is pre-computed.
This option is suitable for mission configurations that involve CubeSats that need to maintain formation with tight tolerances.
The tight tolerances for a given formation can lead to scaling issues and truncation errors when solving low-thrust trajectory optimization problems using the \texttt{`absolute'} option.
Given that the \texttt{orbit\_model} option \texttt{`relative'} is selected in Listing~\ref{lst:vt}, the initial orbit state is defined as a position relative to the initial state of the reference orbit.

\begin{listing}
\begin{minted}[]{python}
cs = CubesatSpec(
    name='optics',
    groundstations=groundstations,
    orbit_model='relative',
    attitude=None,
    dry_mass=1.3, # kg
    initial_orbit_state=np.array([0.02, 0.001, -0.003, 0., 0., 0.]), # km
    specific_impulse=47., # s
)
\end{minted}
\caption{Creating a specification for a CubeSat configuration. Note the use of a list of ground station specifications from Listing~\ref{lst:gs} to specify that the CubeSat is in communication with those ground stations.}
\label{lst:cs}
\end{listing}

Specifications can also control whether or not a particular discipline is modeled.
In Listing~\ref{lst:cs}, the \texttt{attitude} parameter is set to \texttt{None}, indicating that the attitude dynamics of the CubeSat are not modeled for this problem.
Note that if the \texttt{groundstations} object is an empty dictionary, the communication discipline is not modeled.

\subsubsection{Virtual Telescope Specification}\label{sec:vt}

\toolbox{} provides a specification for a specific type of swarm model: a virtual telescope.
A virtual telescope is a telescope whose components are distributed across different locations rather than a single device.
In the case of space-based virtual telescopes, rather than a single spacecraft bus, two spacecraft house the optics and detector.
The advantage of using a virtual telescope in space is the reduction in mass and volume compared to using a single spacecraft bus, reducing launch costs.
The spacecraft need to maintain formation with tight tolerances, however, in order for the mission to be successful.
Listing~\ref{lst:vt} shows a specification for a virtual telescope mission.
The ground station and CubeSat specifications (Sections~\ref{sec:gs}, \ref{sec:cs}, respectively) are among the parameters necessary to specify the system model.
Note that the mission is fully specified without a \toolbox{} user ever writing CSDL code to define the physical behavior of the discipline models.

\begin{listing}
\begin{minted}[]{python}
vt = VirtualTelescopeSpec(
    # dictionary of CubeSat specifications for each CubeSat in the virtual telescope
    optics_cubesat=cubesats['optics'],
    detector_cubesat=cubesats['detector'],
    groundstations=groundstations,
    telescope_length_m=40.,
    telescope_length_tol_mm=0.15,
    telescope_view_halfangle_tol_arcsec=90.,
    max_separation_all_phases_km=5.,
)
\end{minted}
\caption{Creating a specification for a virtual telescope configuration. Note the use of a list of ground station specifications from Listing~\ref{lst:gs} to specify that the CubeSats are communication with those ground stations. The ground station specifications are required at the virtual telescope level because the virtual telescope discipline contains disciplines for the ground stations as well as the communication between each ground station and each CubeSat.}
\label{lst:vt}
\end{listing}

The following subsections present the various discipline models and examples of analyses performed by running the generated computational models directly.

\subsection{Configuration Models Provided by \toolbox{}}

The \texttt{configurations/} directory contains hierarchically composed models of specialized disciplines.
For example, the \texttt{configurations/} directory contains model definitions that capture the vehicle dynamics of a CubeSat and a model of the telescope configuration for a virtual telescope.
The vehicle dynamics configuration contains other disciplines, such as orbit mechanics, attitude dynamics, and propulsion.
The virtual telescope configuration models the in-flight formation of two CubeSats in Low Earth Orbit (LEO); one CubeSat carries the optics and the other carries the detector.

\section{Discipline Models Provided by \toolbox{}}\label{sec:disciplines}

This section presents a subset of the available models contained in the \texttt{disciplines/} directory (see Section~\ref{sec:tour}).
This section is not meant to be a comprehensive guide of all the models available for composing a system model, but it 
describes the most commonly used models that should allow the reader to get started using \toolbox{}
for simple problems.

\subsection{Orbit}\label{sec:orbit}

The orbit dynamics model~\eqref{eq:gravity} uses the $J_2$, $J_3$, and $J_4$ perturbations to capture the effect of the oblateness of the Earth\footnote{Eagle, C.D., ``Orbital Mechanics with MATLAB" http://www.cdeagle.com/ommatlab/toolbox.pdf [retrieved
Februray 2013].}:

\begin{equation}
\begin{aligned}
  \ddot{\vec{r}}=-\frac{\mu}{r^{3}} \vec{r} &-\frac{3 \mu J_{2} R_{e}^{2}}{2 r^{5}}\left[\left(1-\frac{5 r_{z}^{2}}{r^{2}}\right) \vec{r}+2 r_{z} \hat{z}\right] \\
  &-\frac{5 \mu J_{3} R_{e}^{3}}{2 r^{7}}\left[\left(3 r_{z}-\frac{7 r_{z}^{3}}{r^{2}}\right) \vec{r}\right.\\
  &\quad\quad\quad\quad\quad
  +\left.\left(3 r_{z}-\frac{3 r^{2}}{5 r_{z}}\right) r_{z} \hat{z}\right] \\
  &+\frac{15 \mu J_{4} R_{e}^{4}}{8 r^{7}}\left[\left(1-\frac{14 r_{z}^{2}}{r^{2}}+\frac{21 r_{z}^{4}}{r^{4}}\right) \vec{r}\right.\\
  &\quad\quad\quad\quad\quad
  +\left.\left(4-\frac{28 r_{z}^{2}}{3r^{2}}\right) r_{z} \hat{z}\right] .
\end{aligned}\label{eq:gravity}
\end{equation}

\begin{figure}[H]
\centering
\includegraphics[width=0.75\textwidth]{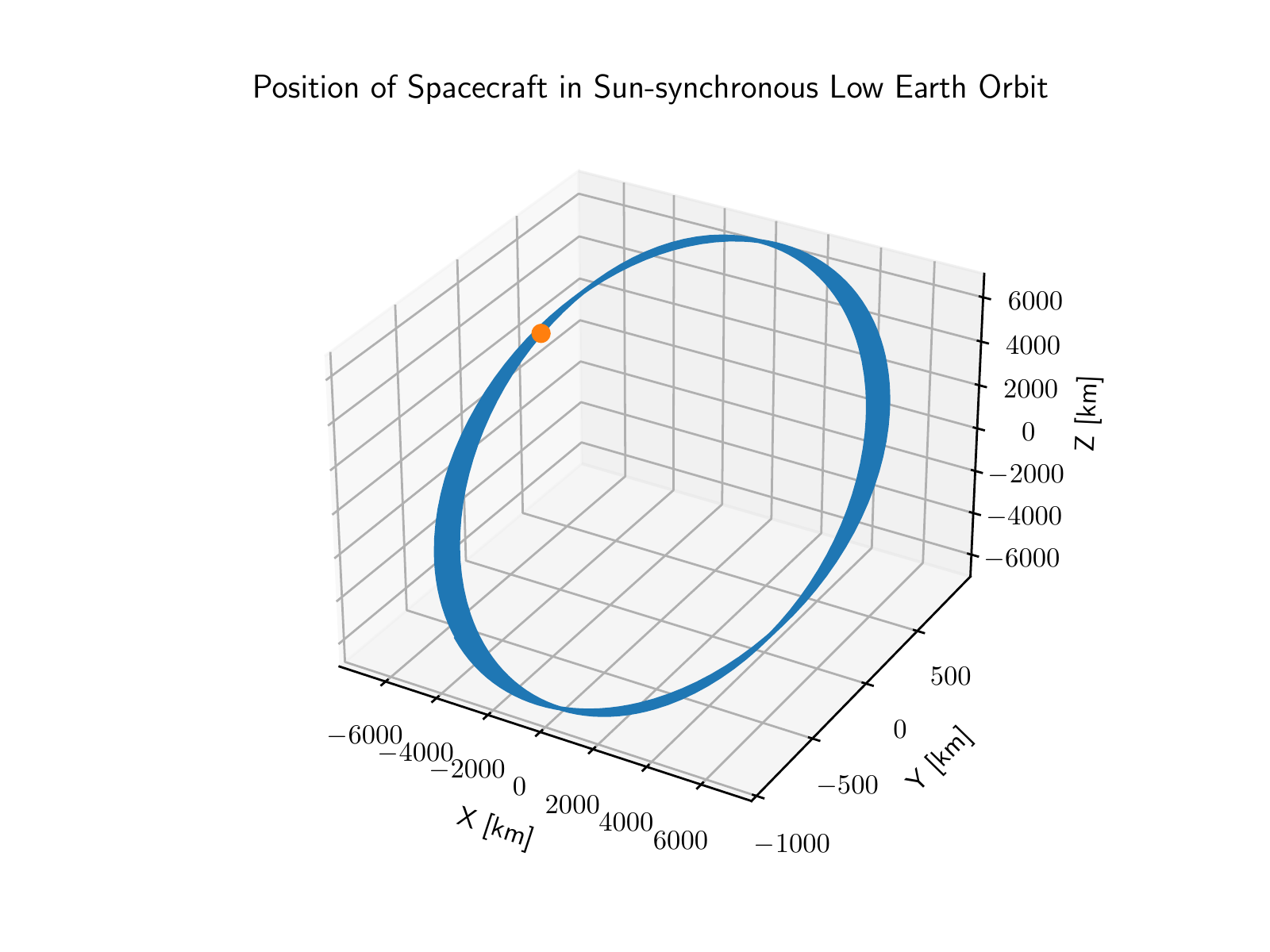}
\caption{Position of a single spacecraft in a sun-synchronous reference orbit taking into account $J_2$, $J_3$, and $J_4$ perturbations over 30 orbits. Coordinates are in the Earth-Centered Inertial frame of reference.}
\label{fig:orbit}
\end{figure}

Figure~\ref{fig:orbit} shows the trajectory of a spacecraft in orbit about the Earth under the influence of gravity, taking into account $J_2$, $J_3$, and $J_4$ perturbations.

\toolbox{} also models the motion of spacecraft relative to a reference orbit, while taking into account $J_2$, $J_3$, and $J_4$ perturbations.
The position of the $i^{\text{th}}$ spacecraft $r^{(i)}$ is given by

\begin{align}
r^{(i)} = r^{(0)}+u^{(i)},
\end{align}

\noindent where $r^{(0)}\in\mathbb{R}^3$ is the reference orbit position and $u^{(i)}\in\mathbb{R}^3$ is the position of the $i^{\text{th}}$ spacecraft relative to the reference orbit.
When used within an optimization framework (e.g. trajectory optimization), constraints that define a spacecraft formation depend on the $u^{(i)}$.
The large difference in orders of magnitude between the constraint values and $r^{(i)}$ leads to truncation errors during optimization even if double precision floating point values are used to compute the trajectory.
These truncation errors prevent the optimizer from finding a solution to an optimization problem.
To avoid truncation errors, we reformulate the equations of motion in terms of the position of the spacecraft relative to the reference orbit, which is given by

\begin{align}\label{eq:relmot}
\ddot{u}^{(i)} = \ddot{r}^{(i)} - \ddot{r}^{(0)}.
\end{align}

The $\ddot{r}^{(i)}$ term is defined in terms of $r^{(0)}$ and $u^{(i)}$, so $\ddot{u}^{(i)}$ can be defined entirely in terms of $r^{(0)}$ and $u^{(i)}$.
Any constraints that depend on $u^{(i)}$ then depend on the design variables (e.g., a term representing thrust added to Equation~\eqref{eq:relmot} in a trajectory optimization problem) in a way that does not introduce truncation error.
Note that the reference orbit $r^{(0)}$ does not depend on design variables and can be precomputed prior to solving an optimization problem to improve performance.

\begin{figure}[H]
\centering
\includegraphics[width=0.75\textwidth]{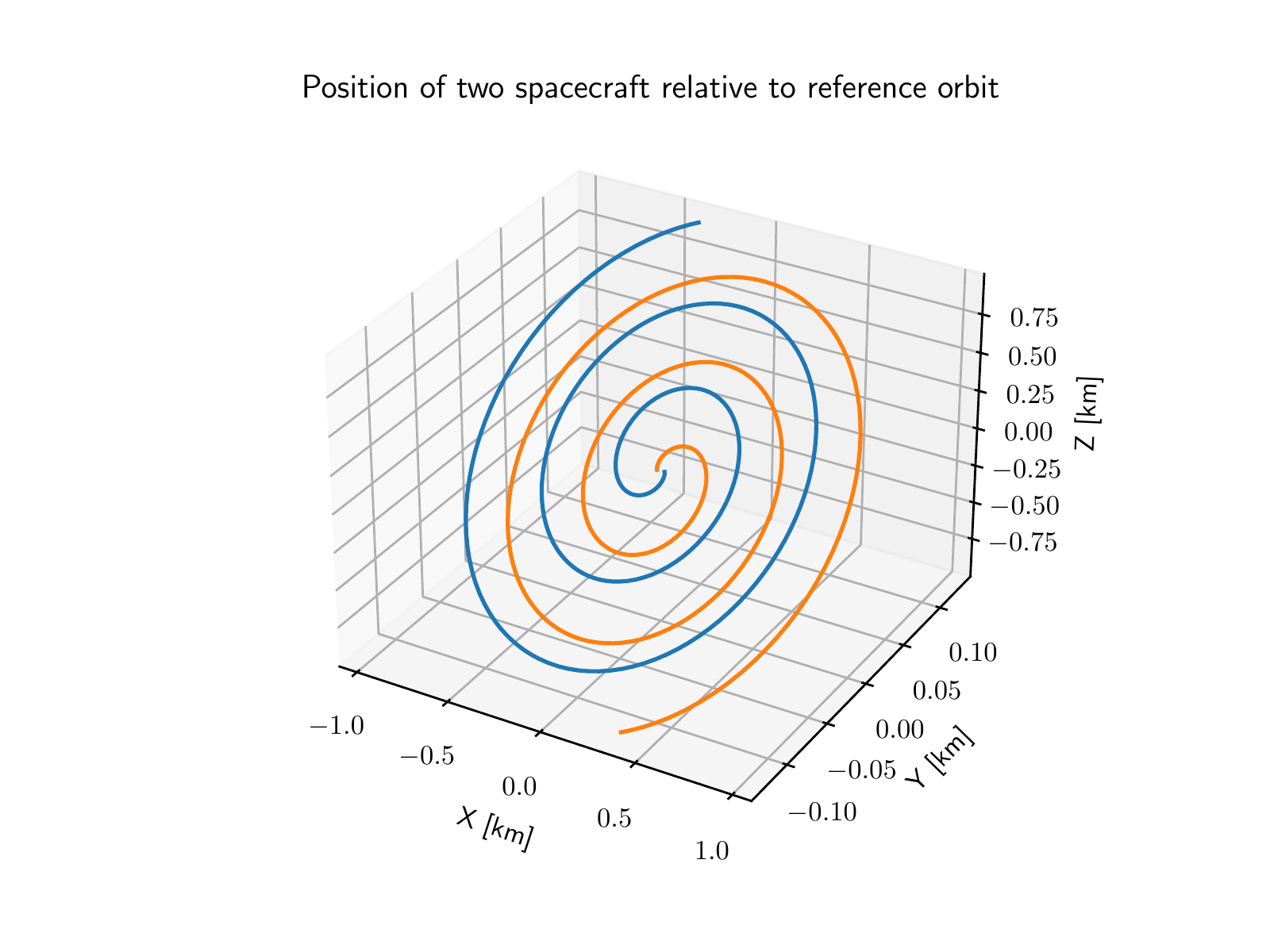}
\caption{Positions of two free flying spacecraft relative to a reference orbit over 3 orbits around the Earth. The initial conditions of the spacecraft are offset from the initial condition for the reference orbit by 20\,m in the $+x$ and $-x$ directions. Coordinates are in the Earth-Centered Inertial frame of reference.}
\label{fig:relative_positions}
\end{figure}

Figure~\ref{fig:relative_positions} shows two trajectories of spacecraft relative to a reference orbit.
The reference frame in Figure~\ref{fig:relative_positions} is aligned with the reference frame in Figure~\ref{fig:orbit}, but the origin of the reference frame in Figure~\ref{fig:relative_positions} is fixed in the location of the spacecraft in Figure~\ref{fig:orbit}.

\subsection{Attitude}\label{sec:attitude}

\begin{figure}[H]
\centering
\includegraphics[width=0.55\textwidth]{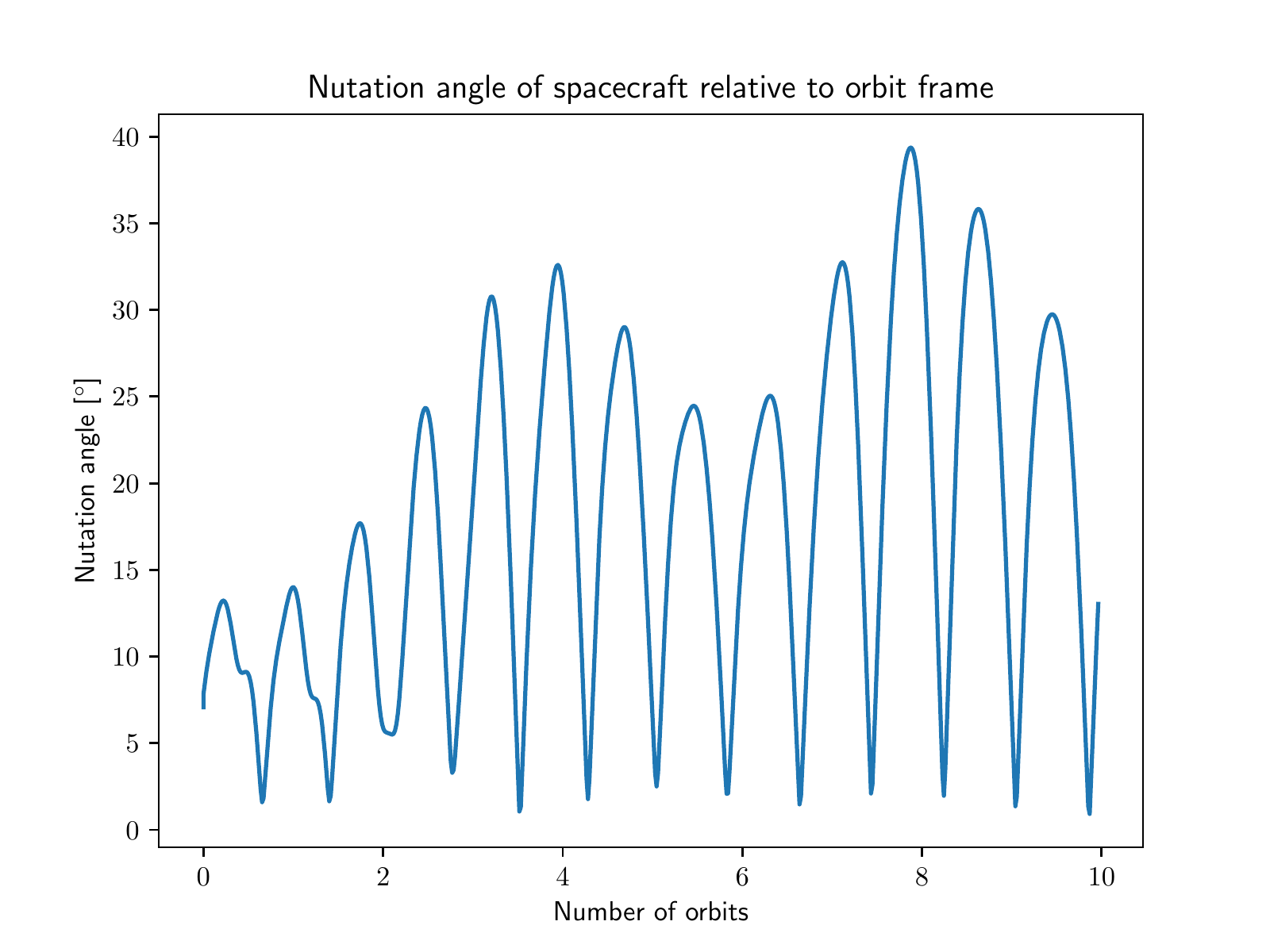}
\caption{Nutation angle of an unsymmetric body in circular orbit under the influence of a gravity gradient assuming Earth is point mass.}
\label{fig:nutation}
\end{figure}

\begin{figure}[H]
\centering
\includegraphics[width=0.75\textwidth]{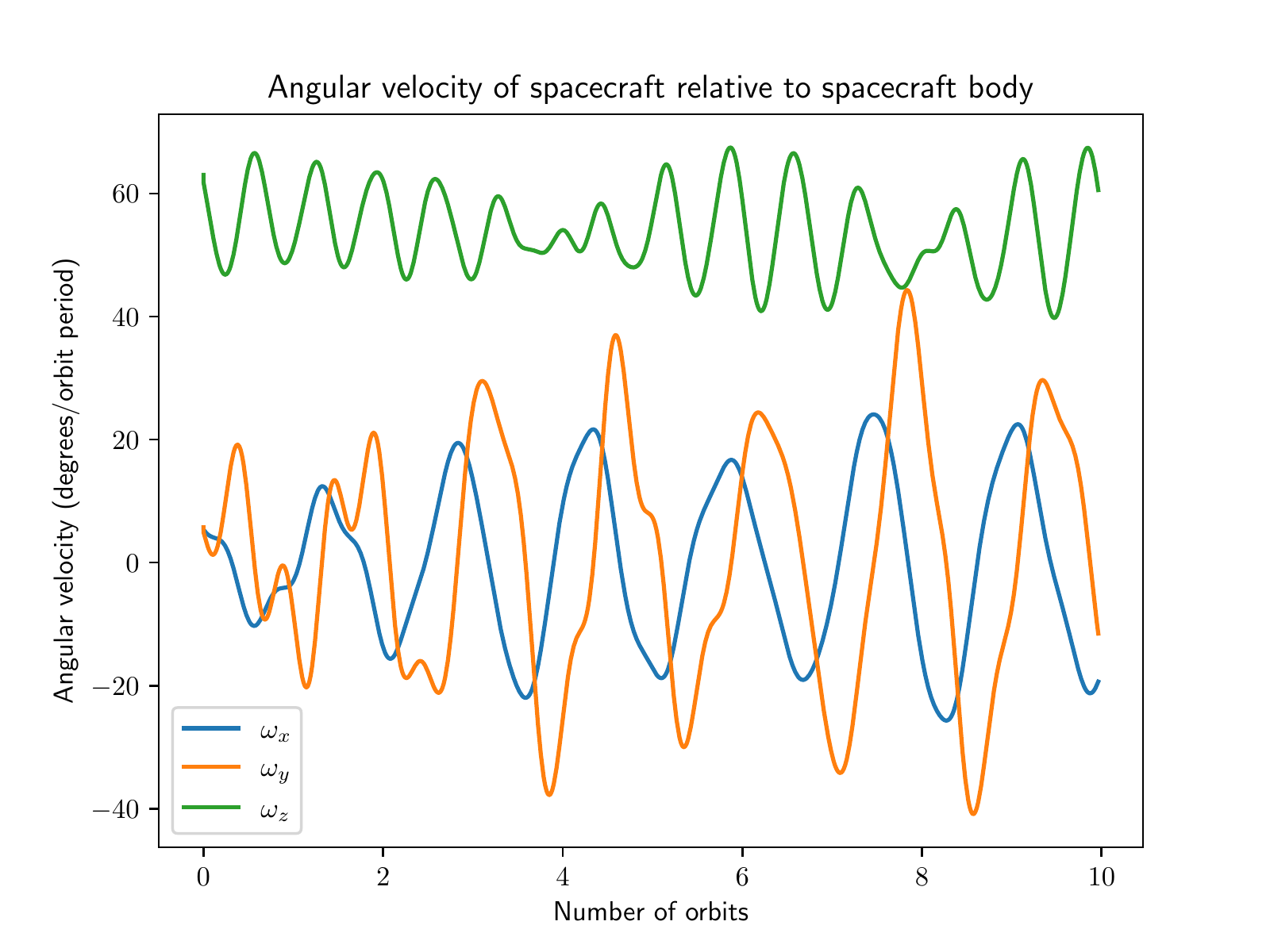}
\caption{Angular velocity of an unsymmetric body in circular orbit under the influence of a gravity gradient assuming Earth is point mass.}
\label{fig:angular_velocity}
\end{figure}

\toolbox{} provides a model of the attitude dynamics of a rigid body in orbit.
Figure~\ref{fig:nutation} shows the nutation angle of a rigid body in circular orbit as a result of the uneven distribution of the force of gravity acting on the body.
The equations of motion of an unsymmetric spacecraft in circular orbit (taken from \cite{kane1983spacecraft}) are given by

\begin{equation}
\begin{aligned}
\dot{\omega}_x &= K_x(\omega_y\omega_z -3\Omega^2C_{21} C_{31})\\
\dot{\omega}_y &= K_y(\omega_z\omega_x -3\Omega^2C_{31} C_{11})\\
\dot{\omega}_z &= K_z(\omega_x\omega_y -3\Omega^2C_{11} C_{21})\\
\dot{C}_{31} &= C_{32}\omega_z - C_{33}\omega_y\\
\dot{C}_{32} &= C_{33}\omega_x - C_{31}\omega_z\\
\dot{C}_{33} &= C_{31}\omega_y - C_{32}\omega_x\\
\dot{C}_{11} &= C_{12}\omega_z - C_{13}\omega_y + \Omega({C}_{13}{C}_{32} - {C}_{12}{C}_{33})\\
\dot{C}_{12} &= C_{13}\omega_z - C_{11}\omega_y + \Omega({C}_{11}{C}_{33} - {C}_{13}{C}_{31})\\
\dot{C}_{13} &= C_{11}\omega_z - C_{12}\omega_y + \Omega({C}_{12}{C}_{31} - {C}_{11}{C}_{32}),
\end{aligned}
\end{equation}

\noindent where $\omega$ is the angular velocity of the spacecraft spinning about its center of mass in the body frame where the axes are aligned with the spacecraft's principal moments of inertia, $\Omega$ is the angular speed of the spacecraft around the central body, and $C$ is the rotation matrix that transforms a vector in the radial-tangential-normal plane to a vector in the body frame.
Additionally, $K_1 = (I_2-I_3)/I_1$, $K_1 = (I_2-I_3)/I_1$, and $K_1 = (I_2-I_3)/I_1$, where $I$ is the vector of moments of inertia about the principal axes of inertia.
Figure~\ref{fig:angular_velocity} shows the angular velocities in the $x$, $y$, and $z$ directions in the body-fixed reference frame.

\subsection{Solar Illumination}\label{sec:illumination}

\begin{figure}[H]
\centering
\includegraphics[width=0.75\textwidth]{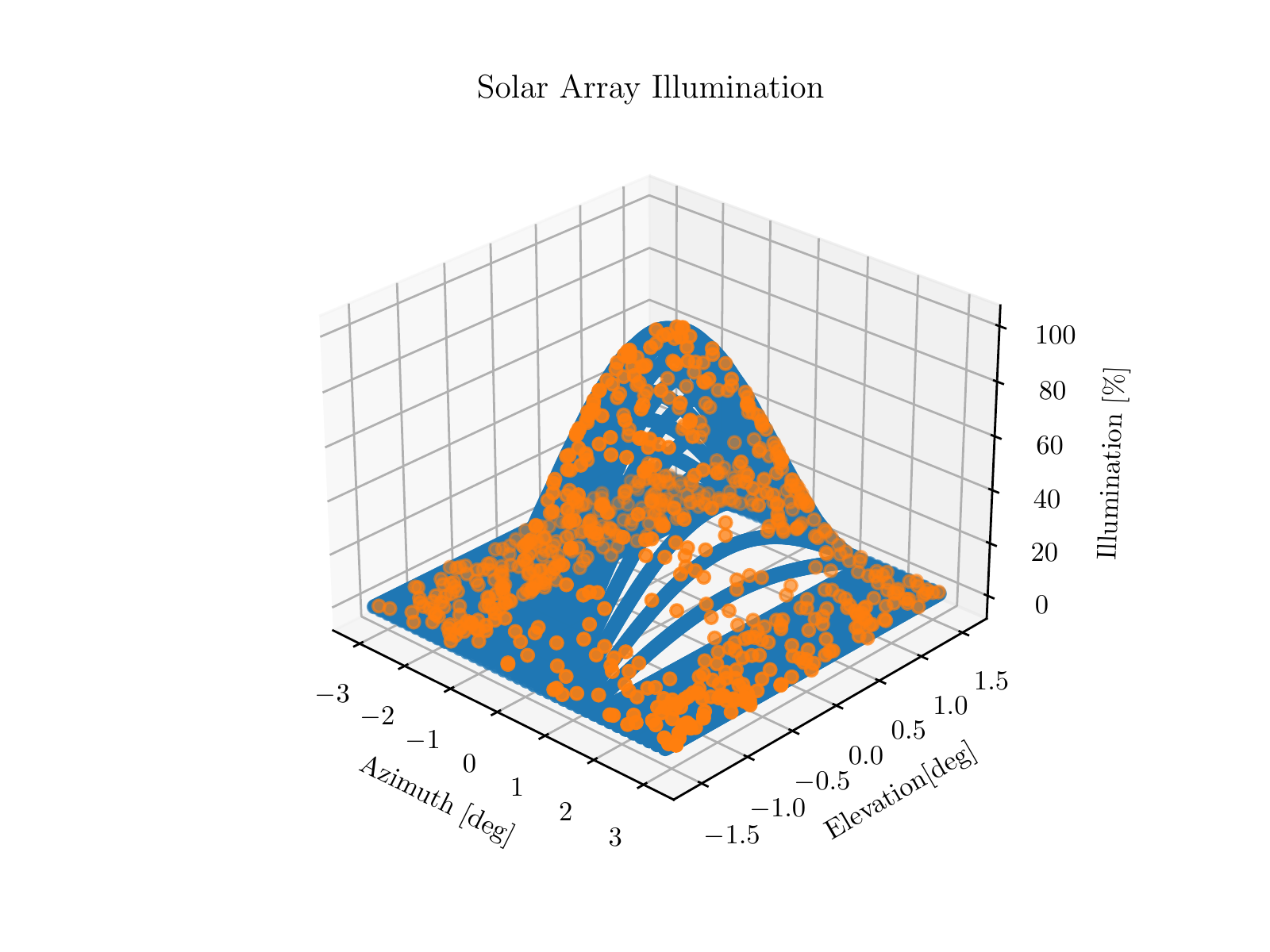}
\caption{Solar illumination as a percent of the area of solar arrays illuminated as a function of azimuth and elevation. Training data (blue) and test data (orange) shown for surrogate model. The training data is generated using a ray tracing algorithm provided by \toolbox{}.}
\label{fig:illumination}
\end{figure}

To model the amount of power generated by solar arrays, the amount of solar illumination on the solar arrays must be known.
The amount of solar illumination is a function of the orientation of the spacecraft relative to the Sun, as well as the shape, orientation, and placement of the solar arrays on the spacecraft.
Ray tracing algorithms are capable of computing the total area of the solar arrays illuminated by the Sun, taking into account shadows cast from one part of a spacecraft onto the solar arrays.
Although ray tracing algorithms can accurately compute the total area of the solar arrays illuminated by the Sun, they are computationally expensive and cannot provide derivatives for use within a gradient-based optimization framework.
\toolbox{} provides utilities for generating training data by applying ray tracing algorithms to CAD files.
These data can then be used to train a surrogate model of the total area of the solar arrays illuminated by the Sun as a continuous and differentiable function of azimuth and elevation angles.
The generated surrogate model is inexpensive to evaluate and also computes derivatives.
Figure~\ref{fig:illumination} shows the training and test data for an example surrogate model using regularized minimal-energy tensor-product splines~\cite{SMT2019} of the total area of the solar arrays illuminated by the Sun as a continuous and differentiable function of the azimuth and elevation angles.

\subsection{Communication}

\begin{figure}[H]
\centering
\includegraphics[width=0.75\textwidth]{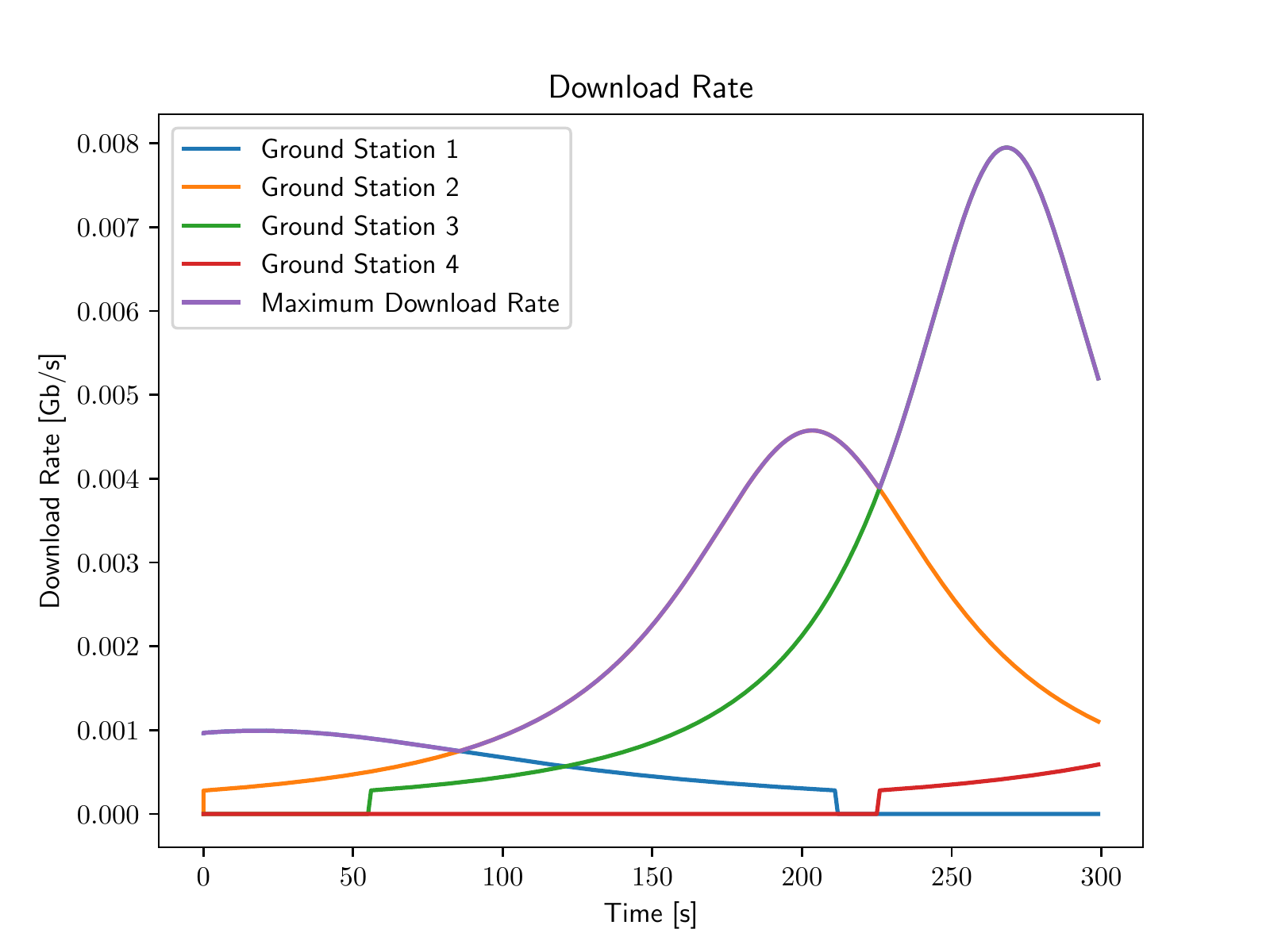}
\caption{Rate of data downloaded from spacecraft for four ground stations, and maximum of data download rate from all ground stations.}
\label{fig:comm}
\end{figure}

The communication discipline model computes the data download rate for all ground stations for each CubeSat for which there is a line of sight with the ground station.
The equation for the data download rate, taken from \cite{larson1992space}, is given by

\begin{align}\label{eq:bitrate}
B_r = \frac{c^2G_rL_l}{16\pi^2f^2kT_s(\text{SNR})}\frac{\eta_pP_{comm}G_t}{S^2}\text{LOS}_c,
\end{align}

\noindent where $c$ is the speed of light, $G_r$ is the receiver gain, $G_t$ is the transmitter gain, $L_l$ is the line loss factor, $f$ is the transmission frequency, $k$ is the Boltzmann constant, $T_s$ is the system noise temperature, $\text{SNR}$ is the signal-to-noise ratio, $\eta_p$ is transmitter amplification efficiency, $P_{comm}$ is the power used by the spacecraft to transmit data with the ground, $S$ is the distance between the ground station and the spacecraft, and $\text{LOS}_c$ indicates that the ground station has a line of sight to the spacecraft.
In the example shown in Figure~\ref{fig:comm}, the transmitter antenna is assumed to have omnidirectional capability, so at any point in flight, the ground station with the highest download rate is assumed to be receiving data.
To compute the download rate for a given satellite, the maximum download rate from all the ground stations is selected.
The total data downloaded is computed by integrating the maximum download rate of all the ground stations.

\begin{figure}[H]
\centering
\includegraphics[width=0.75\textwidth]{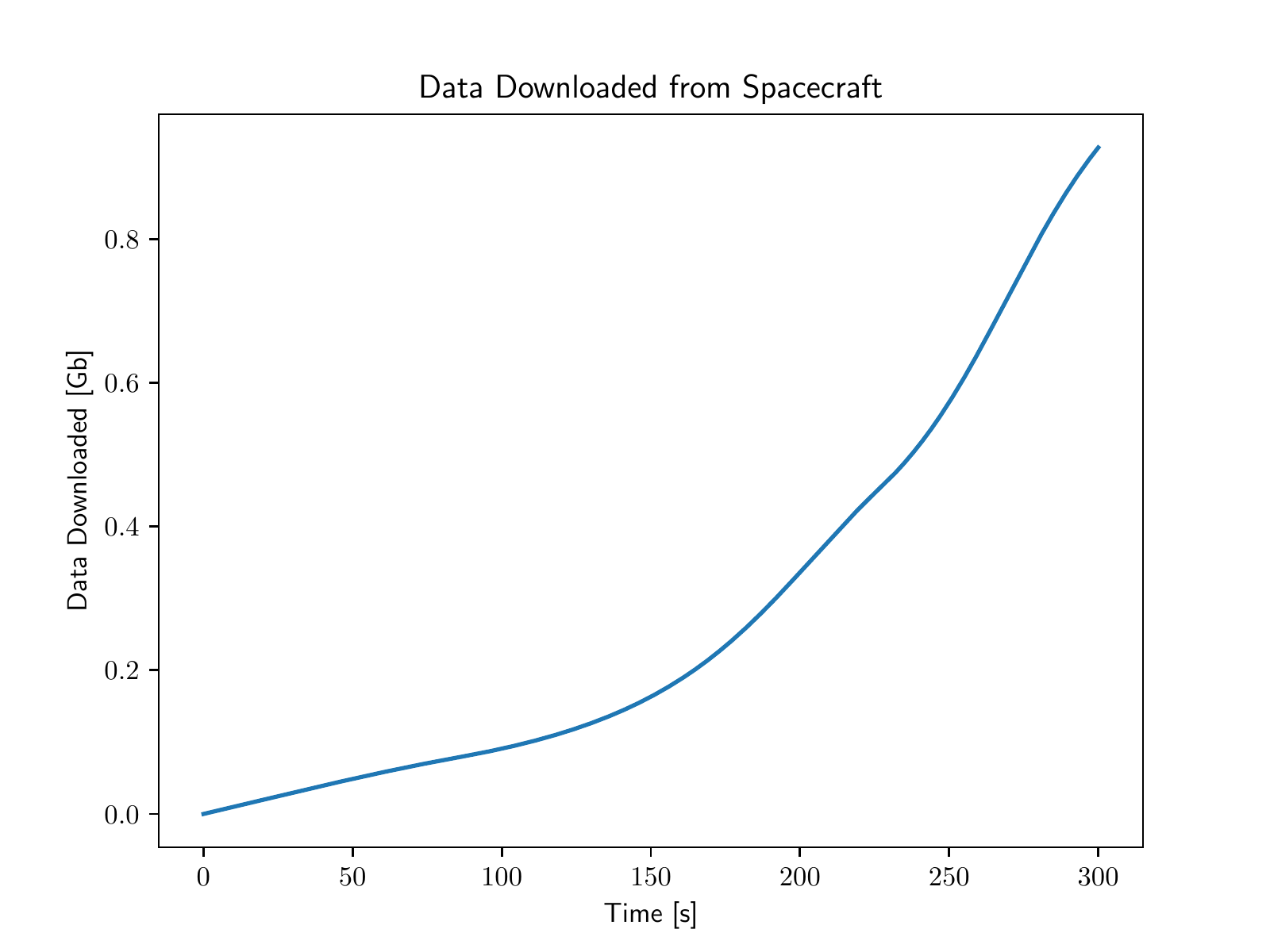}
\caption{Total data downloaded from spacecraft across four ground stations.}
\label{fig:data}
\end{figure}

\section{Trajectory Optimization}

This section presents results from the solution to a trajectory optimization problem for a virtual telescope mission inspired by the VIrtual Super-resolution Optics with Reconfigurable Swarms
(VISORS)~\cite{lightsey2022concept, gundamraj2021preliminary,
koenig2021formation} mission.
The following problem describes the trajectory optimization problem:

\begin{equation}\label{eq:traj}
\begin{aligned}
\min_{T(t)} & \int_{t_0}^{t_f}\sum_{i=1}^2 \sum_{j=1}^3 |\dot{m}_{\text{prop}, i,j}(t)|dt\\
\text{s.t.} &\quad 39.085 ~\text{m}\le \|(u^{(2)}(t) - u^{(1)}(t))\cdot\hat{s}(t)\| I_o(t) \le 40.015~\text{m}, \forall t\in [t_0, t_f] & \text{telescope length}\\
&\quad \arccos\left({\frac{(u^{(2)}(t) - u^{(1)}(t))\cdot\hat{s}(t)}{\|u^{(2)}(t) - u^{(1)}(t)\|}}\right)I_o(t) \le 90~\text{arcsec}, \forall t\in [t_0, t_f] & \text{pointing error}
\end{aligned}
\end{equation}

\noindent where $\dot{m}_{\text{prop}, i,j}$ is the mass flow rate of the propellant used by the $i^\text{th}$ spacecraft in the $j^\text{th}$ direction ($x$, $y$, $z$), and $\hat{s}(t)$ is a unit vector in the direction of the Sun.
The term $I_o(t)$ is an indicator that is 1 when the spacecraft are in the part of the orbit where the telescope makes observations and 0 otherwise.
An additional constraint, $\dot{u}^{(i)}(t) - f(u^{(i)}_0, r^{(0)}(t), T(t)) = 0$, where $u^{(i)}$ are the positions of the spacecraft relative to a reference orbit $r^{(0)}$ (see Section~\ref{sec:orbit}), and $f$ is the differential equation~\eqref{eq:relmot}, may be imposed if using collocation methods.
Otherwise, $f$ can be integrated during each model evaluation.
In the results that follow, the differential equation $f$ is computed explicitly using a Runge-Kutta 4 time integrator.

\begin{figure}[H]
\centering
\includegraphics[width=\textwidth]{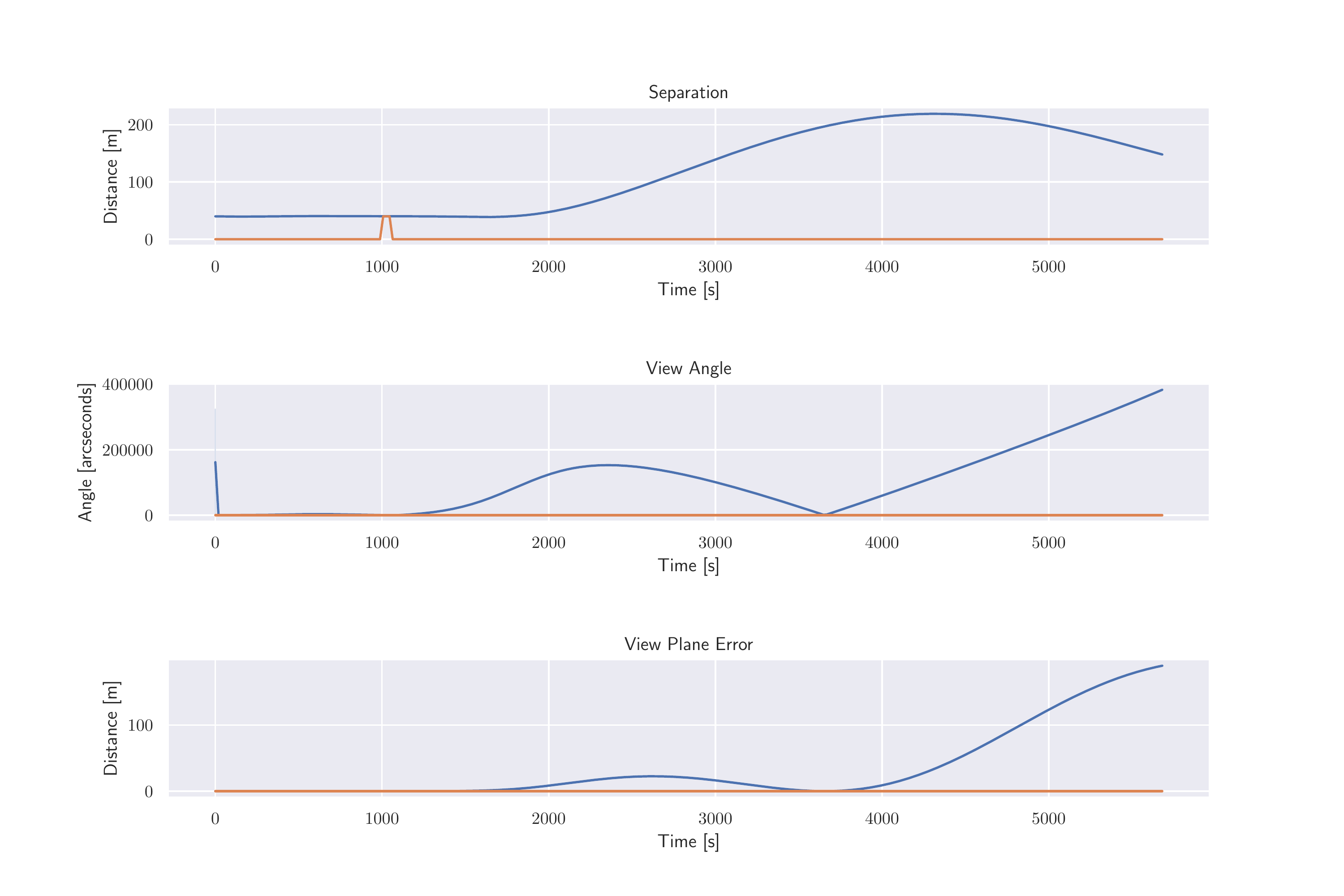}
\caption{Telescope length (spacecraft separation), pointing error, and alignment over approximately one orbit.}
\label{fig:traj}
\end{figure}

\begin{figure}[H]
\centering
\includegraphics[width=\textwidth]{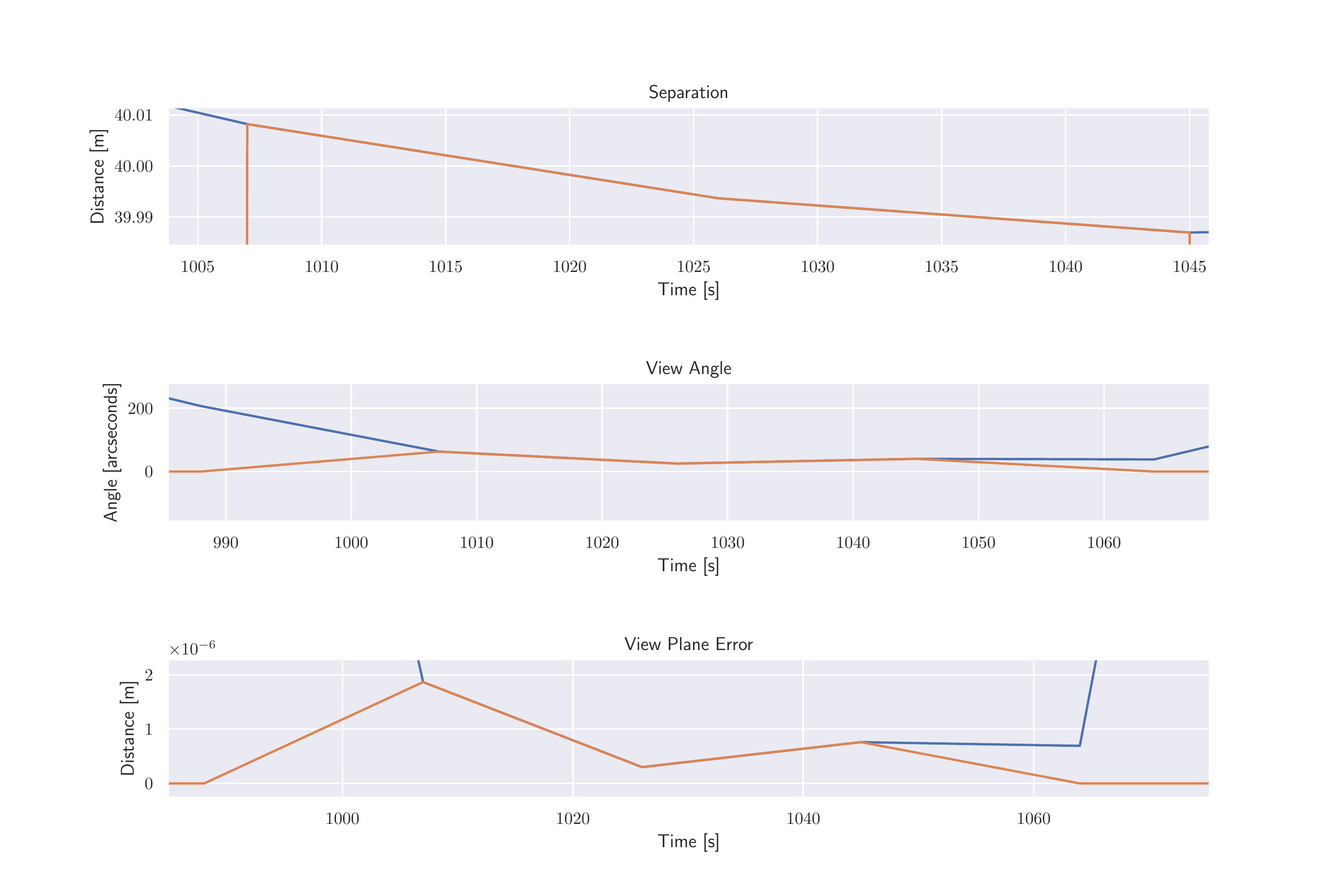}
\caption{Telescope length (spacecraft separation), pointing error, and alignment during observation phase}
\label{fig:traj_zoom}
\end{figure}

\begin{figure}[H]
\centering
\includegraphics[width=\textwidth]{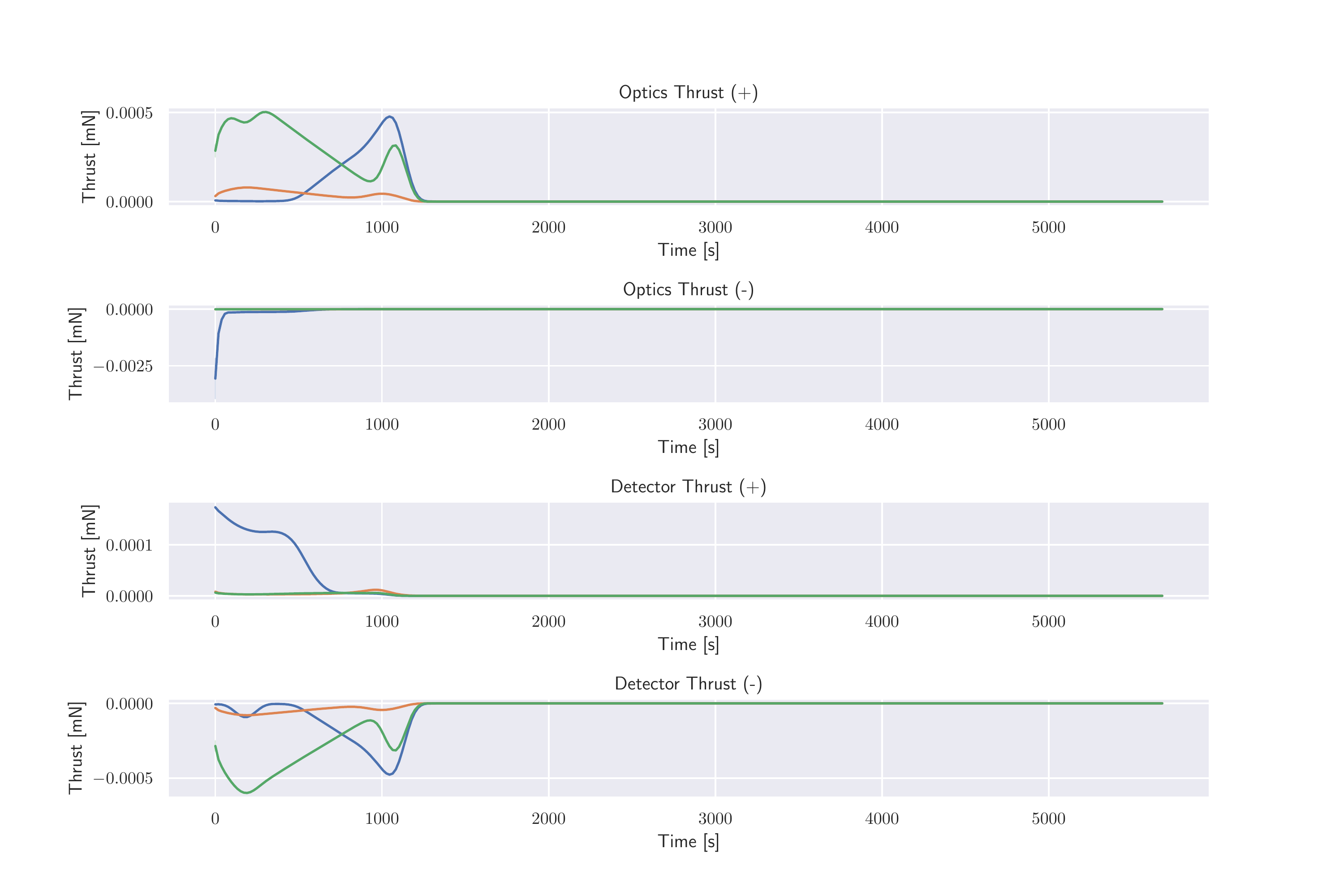}
\caption{Thrust profiles solution $T_i^*(t)$ to \eqref{eq:traj}.}
\label{fig:thrust}
\end{figure}

\begin{figure}[H]
\centering
\includegraphics[width=0.75\textwidth]{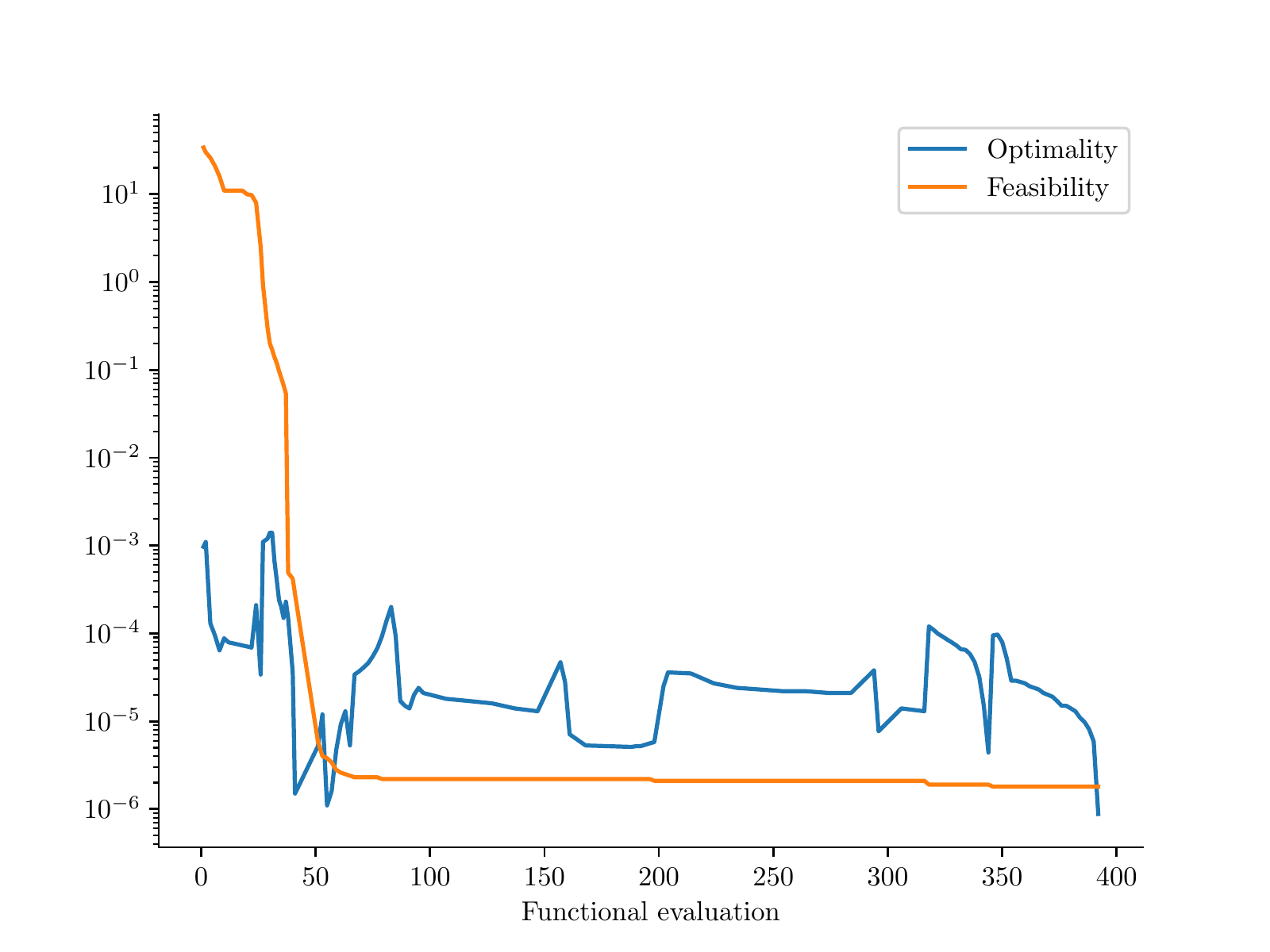}
\caption{SNOPT optimality and feasibility versus iteration number.}
\label{fig:snopt1}
\end{figure}

Figure~\ref{fig:traj} shows the telescope length (spacecraft separation), pointing error, and view plane error (a measure of misalignment between the two spacecraft) over approximately one orbit after the optimization converged.
Figure~\ref{fig:traj_zoom} shows the same information as Figure~\ref{fig:traj} zoomed into the observation phase.
Figure~\ref{fig:thrust} shows the solutions $T_i^*(t)$ to \eqref{eq:traj}.
Figure~\ref{fig:snopt1} shows the feasibility and optimality measures used by SNOPT to determine whether the optimization has converged.
The optimization is considered converged when feasibility is below 1e-3 and optimality is below 1e-6. In Figure~\ref{fig:snopt1}, feasibility falls below 1e-5.

\section{Summary}

This article presented the \toolboxfullname{} (\toolbox{}), a library developed using the Computational System Design Language (CSDL) to explore CubeSat design concepts and conduct trade studies.
CSDL fully automates derivative computation using the adjoint method and separates the definition of the numerical model from the implementation of the computational model.
\toolbox{} provides interfaces for high-level system descriptions, called specifications, for users to fully describe a system model of a CubeSat or CubeSat swarm mission, without the need to write CSDL code directly.
This article provides a tour of the \toolbox{} codebase and presents examples of high-level descriptions (specifications) and discipline models provided by \toolbox{}.
These examples serve to highlight the current set of features of \toolbox{}.

\section{Acknowledgments}

This work was supported in part by the National Science Foundation grant no.~1936557.
The authors thank the undergraduate students at the University of California San Diego who contributed to building the discipline models for \toolbox{}: Andrew Rodriguez, Benjamin Noriega, Joshua Cantrell, Jennifer Nguyen, Tian Nguyen, Etienne Lafayette, Nicholas Lottermoser, Jeffrey Dungan, Brendan Liang, and Ali Kattee.

\printbibliography
\end{document}